\documentclass[10pt,a4paper,reqno,twoside]{article}
\usepackage{hyperref}
\usepackage{amsfonts}
\usepackage{bbding}
\usepackage{multirow}
\usepackage[all]{xy}
\usepackage{amsmath,amscd,amssymb,latexsym,srcltx,indentfirst,titlesec}
\usepackage{hyperref}
\usepackage{amsfonts}
\usepackage{bbding}
\usepackage[all]{xy}
\usepackage{amsmath,amscd,amssymb,latexsym,srcltx,indentfirst,titlesec}
\usepackage{multirow}
\usepackage{booktabs}

\numberwithin{equation}{section}
\newtheorem{theorem}{Theorem}[section]
\newtheorem{lemma}[theorem]{Lemma}

\newtheorem{corollary}[theorem]{Corollary}

\newtheorem{definition}[theorem]{Definition}

\allowdisplaybreaks

\newcommand{\Irr}{\operatorname{Irr}}
\allowdisplaybreaks

\begin{document}

\markboth{Huan Xiong}{Character Graphs Associated with Codegrees}

\title{{\large{\textbf{Finite Groups Whose Character Graphs Associated with Codegrees Have No
Triangles}}}}   %\footnotetext{This work was supported by National 973
%Project (2011CB808003) and NSFC grant (11131001).}}

\author{{\small{
\textbf{ Huan Xiong}} } \\
 {\footnotesize  \textsl{Universit\'e de Strasbourg, CNRS, IRMA UMR 7501, F-67000 Strasbourg, France}} \\ {\footnotesize \textsl{E-mail: xiong@math.unistra.fr}}\\
}

\date{} \maketitle
\noindent\textbf{Abstract.} {\small{Motivated by Problem $164$
proposed by Y. Berkovich and E. Zhmud' in their book \lq\lq Characters of Finite Groups\rq\rq, we
give a characterization of finite groups whose irreducible character
codegrees are prime powers. This is based on a new kind of character
graphs of finite groups associated with codegrees. Such graphs have
close and obvious connections with character  coedgree graphs. For
example, they have the same number of connected components.
By analogy with the work of finite groups whose character graphs (associated with degrees) have no triangles,
we conduct a result of classifying finite groups whose character graphs associated with  codegrees have no triangles in the latter part of this paper.
}}



 \section{Introduction}

 Let $G$ be a nontrivial finite group throughout the paper. Let  $\chi$ be   a character of $G$ and $a(\chi):= \frac{|G:ker\ \chi |}{\chi(1)}$ the so-called character codegree of the character $\chi $. Y. Berkovich and E. Zhmud' proposed a problem in \cite{Berk} which is to study a finite group $G$ such that $a(\chi)$ is a prime power for every $\chi\in\Irr(G)$ (see Problem $164$ of \cite{Berk}, Page $306$). In this paper we  solve this problem by giving a characterization of such groups.

 \begin{theorem} \label{primepower}
 Let $G$ be a finite group. Then all the codegrees of irreducible characters of $G$ are prime powers if and only if $G$ is a $p$-group or a Frobenius group whose order
 has exactly two prime divisors.
  \end{theorem}

 Our method is based on a new kind of character graph $\Gamma(G)$  of a finite group $G$.

  \begin{definition}   The graph $\Gamma(G)$, which is called character graph of $G$ associated with  codegrees, is defined as follows:
  the vertex set is  $\Irr(G)\setminus \{1_G\}$ and there is an edge between $\chi,\psi\in \Irr(G)\setminus \{1_G\}$ if and only if the greatest common divisor $gcd(a(\chi), a(\psi))\neq 1$.
  \end{definition}

 There are several kinds of graphs associated with irreducible characters and  degrees of finite groups. They were received far-reaching attention in the last more than twenty years.  For detailed information on these graphs, we refer to \cite{Zhang1,Zhang2}, and a survey article \cite{Lewis}.

%Denote $ccd(G)=\{a(\chi)\mid \chi\in \Irr(G)\ \mbox{and } \chi\neq 1_G\}$.
 Based on the idea of character degree graphs of finite groups, G. Qian, Y.Wang, and H. Wei define the character codegree graph $\Delta(G)$ of a finite group $G$ in \cite{QianW}.  Specifically, the graph $\Delta(G)$ is defined as follows: the
 vertices of $\Delta(G)$   are the
primes dividing the codegree of some nonprincipal irreducible
character of $G$, and  the vertices $p$ and $q$ are
connected by an edge if and only if there exists a codegree of some
nonprincipal irreducible character of $G$ divisible by $pq$. In \cite{QianW}, the
authors develop some properties of the graph $\Delta(G)$. For example, they show  that
 if $\Delta(G)$ is connected, then its diameter is at
most $3$.  Also, they show  that $\Delta(G)$ has at most $2$
connected components, and that $\Delta(G)$ is not connected if and only
if $G$ is Frobenius or $2$-Frobenius.

Just as there are close connections between the  character graph and the character degree graph of a finite group, the character graph associated with codegrees and the character codegree graph  of a finite group are closely related. Obviously, the two graphs have the same number of connected components. Furthermore, by Corollary $3.2$ in  \cite{Lewis} and
results about $\Delta(G)$ in \cite{QianW}, we have
the following result:

\begin{theorem} \label{notconnect} Let $G$ be a finite group and  let $\Gamma(G)$ be the character graph of $G$ associated with codegrees. Then the following statements hold$:$

$(1)$ $\Gamma(G)$ has at most $2$ connected components$;$

$(2)$ If  $\Gamma(G)$ is connected, then its diameter is at
most $4;$ and

$(3)$ $\Gamma(G)$
is not connected if and only if $G$ is Frobenius or $2$-Frobenius.
\end{theorem}

Among all kinds of shapes of character graphs (or character degree
graphs) of finite groups, having no triangles is probably a
distinguishing one. Finite groups whose character graphs associated
with degrees have no triangles are characterized (see \cite{Wu} for
the solvable case and \cite{Li} for the nonsolvable case). We
conduct an analogous work of investigating finite groups whose
character graphs associated with  codegrees have no triangles in the
latter part of this paper. Indeed, we get a classification of such
groups:

\begin{theorem} \label{class} Let $G$ be a nontrivial finite group whose
character  graph associated with  codegrees has no triangles.  Then $G$ is isomorphic
to one of the following groups$:$

-The cyclic group of order $2$ or $3;$

-The alternating group $A_4;$

-The symmetric group $S_3;$

-The dihedral group $D_{10}$ of order $10;$ and

-The nonabelian group $F_{7,3}$ of order $21$.

Conversely, all these groups' character  graphs associated with
codegrees have no triangles.

\end{theorem}

With the list of finite groups whose
character  graphs associated with  codegrees have no triangles, we get the following obvious result:

\begin{corollary}
Let $G$ be a finite group whose
character  graph associated with  codegrees has no triangles.  Then $G$ is solvable.
\end{corollary}

\section{Preliminaries}

In this section, we present some preliminary results needed
later.

By Lemma $2.1$ in $\cite{QianW}$, we have the following result:

\begin{lemma}\label{reduction}   Let $\chi\in\Irr(G)$.

$(1)$ For any $N \unlhd G$ with $N \leq ker \     \chi$, $\chi$ may
be viewed as an irreducible character of $G/N$. The codegree
$a(\chi)$ of $\chi$ is the same whenever $\chi$ is seen as an
irreducible character of $G$ or $G/N$. Furthermore, $a(\chi)$ is
independent of the choice of such $N$. In particular,  $\Gamma(G/N)$
is a subgraph of $\Gamma(G)$.

$(2)$  If $M$ is a subnormal subgroup of $G$ and $\psi$ is an
irreducible constituent of $\chi_M$, then $a(\psi)$ divides
$a(\chi)$.

\end{lemma}

\begin{lemma}\label{direct}  Let $G=H\times K$, $\chi\in\Irr(H)$, and $\psi\in \Irr(K)$. If $gcd(|H|,\ |K|)=1$,
then $a(\chi\times \psi)=a(\chi)a(\psi)$.
\end{lemma}

\noindent\textit{Proof.}   We first claim $ker (\chi \times
\psi)=ker\ \chi \times ker\ \psi$. Indeed, $ker (\chi \times \psi)=
H_1 \times K_1$ for some $H_1 \leq H$ and $K_1 \leq K$ by Corollary
$8.20$ of \cite{Rose}.

It is easy to see that $H_1 \leq ker\ \chi$ and $K_1 \leq ker\
\psi$. Thus  $ker (\chi \times \psi)\leq ker \ \chi \times ker\
\psi$. Obviously,  $ker (\chi \times \psi)\geq ker \ \chi \times ker
\ \psi$. Thus the claim holds.

Now $a(\chi \times \psi)= \frac{|G:ker (\chi \times \psi)|}{\chi(1)\psi(1)}= \frac{|H:ker \ \chi|}{\chi(1)}\cdot \frac{|K:ker\ \psi|}{\psi(1)} =a(\chi)a(\psi)$.   \hspace*{\fill} $\Box$

 \begin{corollary}\label{nilpotent} Let $G$ be a nilpotent group. Suppose $G=P_1 \times¡¡P_2\times \cdots \times P_s $ where $P_i\in Syl_{p_i}(G)$
 for $i=1,2,\cdots s$. If $\theta_i\in \Irr(P_i)$ where $i=1,2,\cdots s$, then $a(\theta_1\times \theta_2\times \cdots \times \theta_s)=a(\theta_1)a(\theta_2)\cdots a(\theta_s)$.
\end{corollary}

\vspace{2ex}

Theorem A of \cite{QianW} implies the following result:

\begin{lemma} \label{divisors} For any prime $p\mid  |G|$, there exists
$\chi\in\Irr(G)$ such that  $p\mid  a(\chi)$.
\end{lemma}

We also need Theorem $8.17$ of \cite{I}:

\begin{lemma} \label{zero} Let $\chi\in\Irr(G)$ and suppose $p\nmid  (|G|/\chi(1))$ for some prime $p$.
Then $\chi(g)=0$ whenever $p\mid o(g)$.

\end{lemma}

The following result is  well-known (see \cite{Thomp}):

\begin{lemma} \label{Frobenius}
If $G$ is a Frobenius group with
Frobenius kernel $N$ and Frobenius complement $H$, then $N$ is a
nilpotent group. Furthermore, if $H$ has even order then $N$ is
abelian.

\end{lemma}

\section{Finite Groups Whose All Character
Codegrees Are Prime Powers}

In this section, we prove Theorem \ref{primepower} and thus solve
Problem $164$ proposed by Y. Berkovich and E. Zhmud' in \cite{Berk}
by
 characterizing  finite groups all of  whose irreducible
character codegrees are prime powers.

\noindent\textbf{\emph{Proof of Theorem} \ref{primepower}}  By
Theorem \ref{notconnect}, we can easily check that if $G$ is a
$p$-group or a Frobenius group whose order has exactly two prime
divisors, then all the irreducible character codegrees of $G$ are
prime powers.

Conversely, suppose that $G$ is a finite group whose all character
codegrees are prime powers. If $\Gamma(G)$ is connected, then  all
the irreducible character codegrees of $G$ must be powers of a fixed
prime number. Thus, by Lemma \ref{divisors}, $|G|$  has only one
prime divisor. This means that $G$ is a $p$-group.

If $\Gamma(G)$ is not connected, then $\Gamma(G)$ have two connected
components by Theorem \ref{notconnect}. Since all character
codegrees of $G$ are prime powers, it follows that  $|G|$ has
exactly two prime divisors by Lemma \ref{divisors}.  Thus $G$ can
not be a $2$-Frobenius group. By Theorem \ref{notconnect},  $G$ is a
Frobenius group and $|G|$ has exactly two prime divisors.
 \hspace*{\fill}  $\Box$

\vspace{2ex}

Theorem \ref{primepower} has a direct corollary:

\vspace{-1ex}

\begin{corollary}
If $G$ is a finite group whose all character codegrees are
prime powers, then $G$  is solvable.
\end{corollary}

\section{Finite Groups Whose Character
 Graphs Associated with  Codegrees Have No Triangles }

We first investigate the character graphs associated with  codegrees
of finite simple groups.

\begin{lemma}  \label{simple}  Let $G$ be a finite simple group. Then $\Gamma(G)$ is a complete
graph.

In particular,  if $\Gamma(G)$ has no triangles, then $G$ is a
cyclic group of order $2$ or $3$.
\end{lemma}

\noindent\textit{Proof.}  Let  $1_G\neq \theta, \psi \in \Irr(G)$.
Since $G$ is a finite simple group, it follows that $ker \
\theta=1$, which means that $a(\theta)=| G:ker \ \theta|/\theta(1)=|
G |/\theta(1)$. Similarly,  we have $a(\psi)=| G:ker \
\psi|/\psi(1)=| G|/\psi(1)$. Also, ${\theta(1)}^2<| G |$, $
{\psi(1)}^2<| G |$ since
$\sum\limits_{\chi\in\Irr(G)}{{\chi(1)}^2}=| G |$. Therefore
$\theta(1)\psi(1)<| G |$, which means there exists some prime $p$
and some positive integer $k$ such that $p^k\mid   | G |$ and
$p^k\nmid \theta(1)\psi(1)$. Then we see that $p\mid a(\theta)$ and
$p\mid a(\psi)$, which means that $a(\chi)$ and $a(\psi)$ are not
coprime and there is an edge between any two vertices in
$\Gamma(G)$, i.e., $\Gamma(G)$ is a complete graph.

If  $G$ is a finite simple group such that $\Gamma(G)$ has no triangles. By the above argument,  $\Gamma(G)$ is a
complete graph.  Therefore, $G$ contains at most $2$ nonprincipal irreducible
characters, which implies that $G$ contains at most $3$ conjugacy
classes. Thus  $G$ is
a cyclic group of order $2$ or $3$. \hspace*{\fill}
$\Box$

\vspace{2ex}

Now we study  finite abelian groups whose character graphs  associated with  codegrees have no triangles.

\begin{lemma} \label{number}  Let $G=C_n$ be a finite cyclic group of order n.
Then for any positive integer $d\mid n$,  $|\{ \chi\in
Irr(G)\mid a(\chi)=d\}|=\varphi(d)$, where $\varphi$ is the
well-known Euler's totient function such that $\varphi(d)=|\{m\mid 0<m \leq d$ and $gcd(m,d)=1\}|$.

\end{lemma}

\noindent\textit{Proof.} Let $G = \langle g\rangle$  and let
$\varepsilon_n$ be a primitive $n$th root of unity. For any $1\leq
i\leq n$, set $\chi_i(g^k)=\varepsilon_n^{ki}$. Then $\chi_i(1\leq
i\leq n)$ are precisely all the irreducible characters of $G$.
Moreover, there are exactly $\varphi(n)$ irreducible characters
$\chi_i$ of G such that $ker \    \chi_i=1$. Notice that for any
$1\leq i\leq n$, $\chi_i(1)=1$. Thus $a(\chi_i)=| G:ker \
\chi_i|/\chi_i(1)=| G:ker \ \chi_i|$. For any positive integer
$d\mid n$, $a(\chi_i)=d$ is equivalent to either $| ker \
\chi_i|=\frac{n}{d}$  or $ker \     \chi_i=\langle g^d\rangle$. Let
$N_d=\langle g^d\rangle$. Then $ker \ \chi_i=N_d$ if and only if
$\chi_i$ can be seen as a faithful irreducible character of $G/N_d$.
Since the number of faithful irreducible characters of $G/N_d$ is
$\varphi(d)$, we have $|\{ \chi\in\Irr(G)| a(\chi)=d\}|=\varphi(d)$.
\hspace*{\fill} $\Box$

 \begin{lemma} \label{abelian} Let $G$ be an abelian group such that  $\Gamma(G)$ has no triangles. Then $G$ is a cyclic group of order
 $2$ or $3$.
\end{lemma}

\noindent\textit{Proof.} By Corollary \ref{nilpotent} and Lemma
\ref{number}, we may assume that $G$
 is a $p$-group. Furthermore, if $|G|>3$, then $G$ has three distinct nonprincipal irreducible
characters whose codegrees have common divisor $p$ by Lemma
\ref{number} and thus $\Gamma(G)$ has a triangle , a contradiction.
 Thus  $G$ is a cyclic group of order
 $2$ or $3$.
 \hspace*{\fill} $\Box$

 \vspace{2ex}

Next we give reduction results about finite groups whose character graphs associated with  codegrees have no triangles.

 \begin{lemma} \label{reduction2}   Let $G$ be a nontrivial finite group such that  $\Gamma(G)$ has no triangles. If $N$ is a maximal normal subgroup of $G$, then $|G/N|=2$ or $3$, $G'=N$,  and $N$ is the unique maximal normal subgroup of $G$.
\end{lemma}

\noindent\textit{Proof.} If $G'=1$, then the lemma holds by Lemma
\ref{abelian}. So we may assume that $G$ is nonabelian. Note that
$\Gamma(G/N)$ is a subgraph of $\Gamma(G)$. Since $\Gamma(G)$ has no
triangles, it follows that $|G/N|=2$ or $3$ by Lemma \ref{simple}.
Thus $G'\leq N$. But since $\Gamma(G/G')$ also has no triangles, we
have $|G/G'|=2$ or $3$ by Lemma \ref{abelian}. Thus $G'=N$.

If $M$ is another maximal normal subgroup of $G$, then we have $M=G'=N$ by the above argument. This means that $N$ is the unique maximal normal subgroup of $G$.   \hspace*{\fill} $\Box$

\begin{lemma} \label{reduction3}    Let $G$ be a nonabelian finite group
such that $\Gamma(G)$ has no triangles. If $G$ has a normal subgroup
$N$ such that $| G:N|=3 $, then $G$ is isomorphic to the alternating
group $A_4$ or the nonabelian group $F_{7,3}$ of order $21$.
\end{lemma}

\noindent\textit{Proof.} Note that  $N\neq1$ since $G$ is
nonabelian. Let $Irr(G/N)=\{\chi_1=1_G, \chi_2, \chi_3\}$ and $P\in
Syl_3(G)$. By Lemma \ref{reduction2}, we may assume that $G'=N$ and
that $N$ is the unique maximal normal subgroup of $G$.

\vspace{2ex}

\textbf{Step $1$.}  We claim that $3\nmid  | N|$ and
thus $| P|=3$. Otherwise, suppose that $3\mid | N|$.
Then by Lemma \ref{divisors},  there exists $\theta\in\Irr(N)$
such that $3\mid a(\theta)$. Let $\chi$ be an irreducible
constituent of $\theta^G$. Then $3\mid a(\chi)$ by Lemma \ref{reduction}.
Since $3\mid a(\theta)$, we have $\theta\neq1_N$. Thus $\chi_N\neq
1_N$ and $\chi \notin\Irr(G/N)$. Therefore we obtain $\chi, \chi_2$, and $\chi_3$ which are three distinct irreducible characters of $G$ and $3\mid
a(\chi)$, $3\mid a(\chi_2)$, $3\mid a(\chi_3)$. Then it follows that
$\Gamma(G)$ contains a triangle, a contradiction. Thus $3\nmid \
| N|$ and $| P|=3$,  as claimed.

\vspace{2ex}

\textbf{Step $2$.} Notice that $G=NP$, $N\bigcap P=1$. We claim that $G$ is a Frobenius group with
Frobenius kernel $N$ and Frobenius complement $P$. Let $\psi\in
\Irr(G)-\Irr(G/N)$. Then $ker \    \psi\leq N$ since $N$ is the
unique maximal normal subgroup of $G$. By Step $1$, we have
$3\nmid \  | N|$ and hence $3\nmid \  | ker \ \psi|$. Since
$3\mid a(\chi_2)$, $3\mid a(\chi_3)$ and $\Gamma(G)$ contains no
triangles, we conclude that $3\nmid \ a(\psi)=\frac{| G|}{|
ker \    \psi\mid \psi(1)}$. Hence $3\nmid \   \frac{| G|}{
\psi(1)}$ since $3\nmid \  | ker \ \psi|$. Therefore, by Lemma
\ref{zero}, for any $1\neq x\in P$, we obtain $\psi(x)=0$. By the
second orthogonality relation,  we have $| C_G(x)|=\sum\limits_{\chi\in
Irr(G)}{\chi(x)\chi(x^{-1})}=\sum\limits_{\chi\in
Irr(G/N)}{\chi(x)\chi(x^{-1})}=3=| P|.$ But since $| P|=3$, $C_G(x)\geq P$
 and thus $C_G(x)=P$ for any $1\neq x\in P$.  Therefore,  $G$
is a Frobenius group with Frobenius kernel $N$ and Frobenius
complement $P$, as claimed.

\vspace{2ex}

\textbf{Step $3$.} We claim that $N$ is a $q$-group for some prime
$q$. Otherwise, we may assume that $N=Q\times R$ where $Q$ and $R$
are the Sylow $q$- and $r$- subgroups of $N$ respectively by lemma
\ref{Frobenius} (the proof is similar if $| N|$ has more than two
prime divisors). Since $|P|=3$ and $G$ is a Frobenius group by the
Step $2$, either $q \geq 5$ or $r \geq 5$. We may assume $r \geq 5$.
By Step $2$, $G$ is a Frobenius group with Frobenius complement $P$
of order $3$. This implies that
 every  nonprincipal irreducible
character of $N$ has precisely $3$ distinct  conjugates in $G$. Then
we may choose irreducible characters $1_R \neq \xi_i(i=1,2)$ of $R$
such that $\xi_1$ and $\xi_2$ are not conjugate in $G$ since $R$ has
at least $4$ distinct nonprincipal irreducible characters.

Let $\theta\in \Irr(Q)$ such that $q \mid a(\theta)$. Then we have
three irreducible characters of $N:$ $\theta \times 1_R$, $\theta
\times \xi_1$, and $\theta \times \xi_2$. It follows that
$(\theta\times 1_R)^G, ( \theta\times \xi_1)^G$, and $(\theta \times
\xi_2)^G$ are three distinct irreducible characters of $G$ since $G$
is a Frobenius group. By Lemmas \ref{reduction} and \ref{nilpotent},
$q$ divides the codegrees of these three characters, which means
that $\Gamma(G)$ has a triangle, a contradiction. Therefore, $N$ is
a $q$-group for some prime $q$, as claimed.

\vspace{2ex}

\textbf{Step $4$.} We claim that $q\leq 7$. Otherwise, suppose
$q>7$.
 Since $|
Z(N)|\geq q>7 $,  $N$ has at least $8$ conjugacy classes. Hence $N$ has at least $7$ nonprincipal
irreducible characters whose codegrees are divided by $q$. Since $|P|=3$,  there are at least
$3$ nonprincipal irreducible characters $\theta_1,\theta_2,\theta_3$
of $N$ which are not conjugate to each other in $G$ and $q\mid \  |
a(\theta_i)|(1\leq i\leq3)$. Therefore,
${\theta_1}^G,{\theta_2}^G,{\theta_3}^G$ are $3$ distinct
irreducible characters of $G$, and $q\mid \  |
a({\theta_i}^G)|(1\leq i\leq3)$ by Lemma \ref{reduction}, which implies that
$\Gamma(G)$ contains a triangle, a contradiction. Thus $q\leq 7$, as claimed.

\vspace{2ex}

\textbf{Step $5$.}
We claim that if $q=5$ or $7$,
 then $| N| = q$.  Otherwise, suppose that $| N|\geq
q^2$. If $N$ is abelian, then $N$ has at least $q^2>7$ conjugacy
classes. If $N$ is nonabelian, then $| Z(N)|\geq q$ and $N$ has at
least $q-1$ noncentral conjugacy classes. Thus $N$ has at least
$2q-1>7$ conjugacy classes. Therefore, whether $N$ is abelian or
not, $N$ has at least $8$ conjugacy classes and  thus has at least
$7$ nonprincipal irreducible characters. By the same discussion as
in Step $4$, we obtain that $\Gamma(G)$ contains a triangle, a
contradiction. Thus if $q=5$ or $7$, then $| N|= q$, as claimed.

\vspace{2ex}

\textbf{Step $6$.} We claim $q\neq 5$. Since  $G$ is a Frobenius
group with Frobenius kernel $N$ and Frobenius complement $P$ of
order 3, $3\mid(| N|-1)$. Thus, $q\neq5$, as claimed.

\vspace{2ex}

\textbf{Step $7$.} We claim that $|N|= 4$ if $N$ is a $2$-group. We
may assume that $| N|=4^n$ for some positive integer $n$ since
$3\mid(| N|-1)$.

Suppose $n=2$ and $|\Irr(N)|\geq 8$. Using the same discussion
as in Step $4$, we deduce that $\Gamma(G)$
contains a triangle, a contradiction.
We now suppose $n=2$ and  $|\Irr(N)|\leq7$. Then  $N$ is
not abelian. By the character table of all nonisomorphic
nonabelian groups of order $16$ in section $25$ of \cite{James},
there are exactly $3$ distinct nonisomorphic nonabelian groups of
order $16$ who have at most $7$ irreducible characters and they all
have precisely three irreducible characters whose degrees are $2$.
Moreover, two  of them are
faithful and the  other one is not faithful. So if $N$
is one of the $3$ distinct nonisomorphic nonabelian groups of order
$16$, then the exact three irreducible characters  of  degree
$2$ can not be conjugate in $G$, which contradicts to  the fact that $G$
is a Frobenius group with Frobenius complement $P$ of order $3$.
Thus $n\neq 2$. So we may assume
 $n\geq3$.

  Indeed,  we may further assume $|\Irr(N)|\leq7$. Since $N$ is nilpotent,
$N'<N$. Note that $G$ is a Frobenius group with Frobenius complement
$P$ of order $3$. Hence the degrees of irreducible characters of $N$
must be $1,1,1,1,y,y,y$ for some positive integer $y$. Therefore, we
have
 $4+3y^2=4^n$.  It
is easy to see $y=2z$ for some positive integer $z$.  So
$1+3z^2=4^{n-1}$ and hence $(2^{n-1}+1)(2^{n-1}-1)=3z^2$. Since
$(2^{n-1}+1,2^{n-1}-1)=1$, we know $2^{n-1}+1=u^2$ or
$2^{n-1}-1=u^2$ for some positive integer $u$. But since $n\geq3$,
$2^{n-1}-1\equiv 3$ (mod $4$), and $u^2\equiv 0$ or $1$ (mod $4$), it follows that
$2^{n-1}-1\neq u^2$. Hence we have $2^{n-1}+1=u^2$, which means
that $(u+1)(u-1)=2^{n-1}$. Thus $u+1$ and $u-1$ are powers of 2. Note that $(u+1)-(u-1)=2$. Hence it follows that $u-1=2$, $u=3$, and $n=4$. But this
leads to $z^2=21$,  a contradiction. Thus the claim holds.

\textbf{Step $8$.}
 Now we show that $N$ is a $2$-group of order $4$ or a cyclic group of order $7$. If $N$ is the latter case, $G$
is  the nonabelian group $F_{7,3}$ of order 21. Notice that the nonprincipal irreducible character codegrees of
$F_{7,3}$ are $3,3,7,7$. The graph $\Gamma(F_{7,3})$ does
contain no triangles.

So we may assume that $N$ is a $2$-group  of order $4$.
 Then $| G|=12$. By the character
table of all $3$ nonisomorphic nonabelian groups of order $12$ in
\cite{James},   the nonprincipal irreducible character codegrees of
$A_4$ are $3,3,4$ and $A_4$ is the only nonabelian group of order
$12$  whose character graph  associated with codegrees contains no
triangles. Thus the proof of the lemma is finished.
 \hspace*{\fill} $\Box$

 \begin{lemma}  \label{reduction4}   Let $G$ be a nonabelian finite group
such that $\Gamma(G)$ has no triangles. If $G$ has a normal subgroup
$N$ such that $| G:N|=2$, then $G$ is isomorphic to the symmetric
group $S_3$ or the dihedral group $D_{10}$ of order $10$.
\end{lemma}

\noindent\textit{Proof.}
 Note that  $N\neq1$ since $G$ is nonabelian.
Let $Irr(G/N)=\{\chi_1=1_G, \chi_2\}$ and $P\in Syl_2(G)$. Then $2\mid a(\chi_2)$. By Lemma \ref{reduction2}, we may assume that $G'=N$ and that $N$ is the unique maximal normal subgroup of $G$.

\vspace{2ex}

\textbf{Step $1$.}  We claim that for any $\chi\in
\Irr(G)-\Irr(G/N)$, $\chi(1)$ is  even. Otherwise, suppose that
there is a  $\psi_1\in \Irr(G)-\Irr(G/N)$ such that  $\psi_1(1)$ is
an odd integer. Since $2\mid \  | G|$ and
$\sum\limits_{\chi\in\Irr(G)-\Irr(G/N)}{{\chi(1)}^2}=| G |-2$, there
is another $\psi_2\in\Irr(G)-\Irr(G/N)$ such that $\psi_2(1)$ is an
odd integer. Note that  $ker \    \psi_i\leq N(1\leq i\leq2)$ since
$N$ is the unique maximal normal subgroup of $G$ and that  $|
G:N|=2$. Therefore, we have $2\mid \  \frac{| G|}{| ker \ \psi_i|}$.
But since $2\nmid \ \psi_i(1)$, it follows that $2\mid\frac{| G|} {|
ker \ \psi_i|\psi_i(1)}$ and  $2\mid a(\psi_i)(1\leq i\leq2)$. This
implies that $\Gamma(G)$ contains a triangle consisting of $\psi_1,
\psi_2$ and $\chi_2$, a contradiction. Thus the claim holds.

\vspace{2ex}

\textbf{Step 2.} We claim that $G$ is a Frobenius group with
Frobenius kernel $N$ and Frobenius complement $P$. By the claim of
Step $1$,   $$| G |=\sum\limits_{\chi\in
Irr(G)}{{\chi(1)}^2}=\sum\limits_{\chi\in
Irr(G)-Irr(G/N)}{{\chi(1)}^2}+2\equiv2(mod\ 4).$$  This means that
$| G |/2$ is an odd integer. Hence $| P|=2$ and $2\nmid \ \frac{|
G|}{ \psi(1)}$ for any $\psi\in \Irr(G)-\Irr(G/N)$. By Lemma
\ref{zero}, for any $1\neq x\in P$, $\psi\in Irr(G)-Irr(G/N)$, we
have $\psi(x)=0$. Therefore, by the second orthogonality relation,
$| C_G(x)|=\sum\limits_{\chi\in
Irr(G)}{\chi(x)\chi(x^{-1})}=\sum\limits_{\chi\in
Irr(G/N)}{\chi(x)\chi(x^{-1})}=2=| P|.$  This means that $C_G(x)=P$
for any $1\neq x\in P$. Thus $G$ is a Frobenius group with Frobenius
kernel $N$ and Frobenius complement $P$, as claimed.

\vspace{2ex}

\textbf{Step 3.} By Lemma \ref{Frobenius}, $N$ is abelian.  Using
the similar argument as in Steps $3,\ 4$, and $5$ in Lemma
\ref{reduction3}, we deduce that $|N|\leq 5$ and thus $N$ is a
cyclic group of order $3$ or $5$. Therefore $G$  is isomorphic to
the symmetric group $S_3$ or the dihedral group $D_{10}$ of order
$10$. \hspace*{\fill} $\Box$

\vspace{2ex}

Now we can prove Theorem \ref{class}:

\noindent\textbf{\emph{Proof of Theorem} \ref{class}}  Let $G$ be a
group appeared in the theorem. By Lemmas \ref{abelian},
\ref{reduction2}, \ref{reduction3}, and \ref{reduction4}, it
suffices to show that $\Gamma(G)$ contains no triangles. This is
obvious if $G$ is a cyclic group of order  $2$ or $3$. For other
cases, it is shown by the following table:

\begin{center}
\begin{tabular}{|l|l|l|l|l|} \hline
Group               &  Nonprincipal irreducible character codegrees
\\    \hline

$S_3$             &   $2,\ 3$ \\  \hline
$D_{10}$        &    $2,\ 5,\ 5$ \\  \hline
$A_4$            &   $3,\ 3,\ 4$   \\  \hline
$F_{7,3}$      &    $3,\ 3,\ 7,\ 7$  \\  

\hline
\end{tabular}
\end{center}
 \hspace*{\fill} $\Box$

\subsection*{Acknowledgement} We really appreciate the encouragement of Prof.
Jiping Zhang.

\end{document}